# ESTIMATES FOR WAVELET COEFFICIENTS ON SOME CLASSES OF FUNCTIONS

**V. F. Babenko**[1,2] **and S. A. Spector**[2]

Let $\psi_m^D$ be orthogonal Daubechies wavelets that have $m$ zero moments and let

$$W_{2,p}^k = \left\{ f \in L_2(\mathbb{R}): \left\| (i\omega)^k \hat{f}(\omega) \right\|_p \leq 1 \right\}, \quad k \in \mathbb{N}.$$

We prove that

$$\lim_{m \to \infty} \sup \left\{ \frac{|(\psi_m^D, f)|}{\left\| (\psi_m^D)^\wedge \right\|_q} : f \in W_{2,p'}^k \right\} = \frac{\dfrac{(2\pi)^{1/p - 1/2}}{\pi^k} \left( \dfrac{1 - 2^{1-pk}}{pk - 1} \right)^{1/p}}{(2\pi)^{1/q - 1/2}}.$$

Let $L_p = L_p(\mathbb{R})$, $1 \leq p \leq \infty$, be the space of measurable functions $f: \mathbb{R} \to \mathbb{C}$ with finite norm $\|f\|_p$, where

$$\|f\|_p = \|f\|_{L_p(\mathbb{R})} = \left( \int_{\mathbb{R}} |f(x)|^p dx \right)^{1/p} \quad \text{if} \quad p < \infty$$

and

$$\|f\|_\infty = \|f\|_{L_\infty(\mathbb{R})} = \operatorname{vrai\,sup}_{x \in \mathbb{R}} |f(x)|.$$

For $f \in L_p(\mathbb{R})$ and $g \in L_q(\mathbb{R})$, where $p, q \in [1; \infty]$, $\dfrac{1}{p} + \dfrac{1}{q} = 1$, we set

$$(f, g) = \int_{\mathbb{R}} f(x) \overline{g}(x) dx.$$

We consider the following classes of functions $f \in L_2(\mathbb{R})$: For $k \in \mathbb{N}$ and $p \in (1, \infty)$, we set

$$W_{2,p}^k = \left\{ f \in L_2(\mathbb{R}): \left\| (i\omega)^k \hat{f}(\omega) \right\|_p \leq 1 \right\},$$

where

---

[1] Dnepropetrovsk National University, Dnepropetrovsk.

[2] Institute of Applied Mathematics and Mechanics, Ukrainian National Academy of Sciences, Donetsk.

$$\hat{f}(\omega) \;=\; \frac{1}{\sqrt{2\pi}} \int_{\mathbb{R}} f(x) e^{-i\omega x} dx$$

is the Fourier transform of a function $f$. For $p = 2$, we obtain the standard Sobolev classes

$$W_{2,2}^k \;=\; \left\{ f \in L_2(\mathbb{R}) : \left\| f^{(k)} \right\|_2 \le 1 \right\}.$$

For a function $\psi(t) \in L_2(\mathbb{R})$ and numbers $j, k \in \mathbb{Z}$, we set

$$\psi_{j,k}(t) \;=\; 2^{j/2} \psi(2^j t - k).$$

If the system of functions $\{\psi_{j,k}\}_{j,k \in \mathbb{Z}}$ forms an orthonormal basis of the space $L_2(\mathbb{R})$, i.e., any function $f \in L_2(\mathbb{R})$ can be represented in the form of the sum of the series

$$f(t) \;=\; \sum_{i \in \mathbb{Z}} \sum_{j \in \mathbb{Z}} \left( \psi_{j,v}, f \right) \psi_{j,v}(t) \tag{1}$$

convergent in $L_2(\mathbb{R})$, then the function $\psi(t)$ is called an orthogonal wavelet.

Many applications of orthogonal wavelets are based on the investigation of the value of the wavelet coefficients in representations of the form (1), depending both on properties of the wavelet $\psi(t)$ and on the smoothness of the function $f$.

Assume that the wavelet $\psi(t)$ has $k$ zero moments or, which is equivalent, that $\hat{\psi}(\omega)$ has a zero of multiplicity $k$ at zero. We define a function $_k\psi(t)$ by the relation

$$\left( _k\psi(t) \right)^{\wedge}(\omega) \;=\; (i\omega)^{-k} \hat{\psi}(\omega).$$

We set

$$C_{\kappa;p,q}(\psi) \;=\; \sup_{f \in W_{2,p'}^k} \frac{|(\psi f)|}{\|\hat{\psi}\|_q},$$

where

$$p' \;=\; \frac{p}{p-1}.$$

It is easy to see that $C_{\kappa;p,q}(\psi)$ can be represented in the form

$$C_{\kappa;p,q}(\psi) \;=\; \frac{\left\| (_k\psi)^{\wedge} \right\|_p}{\|\hat{\psi}\|_q}. \tag{2}$$

Note that $C_{\kappa;p,q}(\psi)$ is the exact constant in the inequality

$$\left|(\psi_{j,\nu}, f)\right| \leq C_{\kappa;p,q}(\psi) 2^{-j\left(k-\frac{1}{p}+\frac{1}{q}\right)} \|\hat{\psi}\|_q \|\hat{f}(\omega)(i\omega)^k\|_{p'}.$$

Let $m \in \mathbb{N}$. The trigonometric polynomials

$$H_m(\omega) = 2^{-1/2} \sum_{l=0}^{2m-1} h_m(l) e^{il\omega}, \quad h_m(l) \in \mathbb{R},$$

which satisfy the equalities

$$|H_m(\omega)|^2 = \left(\cos^2 \frac{\omega}{2}\right)^m P_{m-1}\left(\sin^2 \frac{\omega}{2}\right),$$

where

$$P_{m-1}(x) = \sum_{k=0}^{m-1} \binom{m-1+k}{k} x^k,$$

are called Daubechies filters (see, e.g., [1], Sec. 16).

A function $\varphi_m^D$ whose Fourier transform has the form

$$\left(\varphi_m^D\right)^{\wedge}(\omega) = \frac{1}{\sqrt{2\pi}} \prod_{l=1}^{\infty} H_m(\omega 2^{-l})$$

is the orthogonal scaling function. A function whose Fourier transform has the form

$$\left(\psi_m^D\right)^{\wedge}(\omega) = e^{-\frac{i\omega}{2}} \overline{H_m\left(\frac{\omega}{2}+\pi\right)} \left(\varphi_m^D\right)^{\wedge}\left(\frac{\omega}{2}\right)$$

is called an orthogonal Daubechies wavelet $\psi_m^D$.

The wavelet $\psi_m^D$ possesses the following properties (see [2], Chap. 6; [1], Sec. 16):

(i) $\operatorname{supp} \psi_m^D = [-(m-1), m]$;

(ii) $\psi_m^D$ has $m$ zero moments;

(ii) there exists $\lambda > 0$ such that $\psi_m^D \in C^{\lambda m}$, where

$$C^\alpha = \left\{ f : \int_{\mathbb{R}} \hat{f}(\omega)(1+|\omega|)^\alpha d\omega < \infty \right\}, \quad \alpha > 0. \tag{3}$$

Furthermore (see, e.g., [3], Sec. 5.5),

$$\left|(\psi_m^D)^\wedge(\omega)\right|^2 = \left|H_m\left(\frac{\omega}{2}+\pi\right)\right|^2 \left|(\varphi_m^D)^\wedge\left(\frac{\omega}{2}\right)\right|^2 = \frac{1}{2\pi}\left|H_m\left(\frac{\omega}{2}+\pi\right)\right|^2 \prod_{l=1}^{\infty}\left|H_m(2^{-l-1}\omega)\right|^2, \quad (4)$$

and, moreover,

$$|H_m(\omega)|^2 = \left(1 - c_m \int_0^\omega \sin^{2m-1} u\, du\right). \quad (5)$$

The following theorem is the main result of the present paper:

**Theorem 1.** *Let $k \geq 0$ be a fixed integer. Then*

$$\lim_{m\to\infty} C_{k;p,q}(\psi_m^D) = \frac{(2\pi)^{1/p-1/q}}{\pi^k}\left(\frac{1-2^{1-pk}}{pk-1}\right)^{1/p}.$$

*For $p = q = 2$, this theorem was proved in [4].*

To prove this theorem, we need the lemma presented below, which was also proved in [4] for $p = 2$. Let

$$\hat\Psi = \frac{1}{\sqrt{2\pi}}\left(\chi_{[-2\pi,-\pi]} + \chi_{[\pi,2\pi]}\right),$$

where $\chi_I$ is the characteristic function of the interval $I$.

**Lemma 1.** *Suppose that $1 < p < \infty$, $k \geq 0$ is an integer, $(\psi_n)$ is a sequence of functions with compact support, and the following conditions are satisfied:*

  (i) *for a certain $\varepsilon$ independent of $n$, $0 < \varepsilon < \pi$, one has*

$$\int_{|\omega|<\varepsilon} |\omega|^{-pk}\left|(\psi_n)^\wedge(\omega)\right|^p d\omega \to 0 \quad as \quad n \to \infty;$$

  (ii) $\left\|(\psi_n)^\wedge - \hat\Psi\right\|_p \to 0$ *as $n \to \infty$.*

*Then*

$$\lim_{n\to\infty}\left\|(_k\psi_n)^\wedge\right\|_p = \frac{(2\pi)^{1/p-1/2}}{\pi^k}\left(\frac{1-2^{1-pk}}{pk-1}\right)^{1/p}. \quad (6)$$

***Proof.*** We represent $\left\|({}_k\psi_n)^\wedge\right\|_p$ in the form

$$\left\|({}_k\psi_n)^\wedge\right\|_p = \left(\int_{\mathbb{R}} \left|(i\omega)^{-k}\right|^p \left|(\psi_n)^\wedge(\omega)\right|^p d\omega\right)^{1/p} = \left(\int_{\mathbb{R}} |\omega|^{-pk} \left|(\psi_n)^\wedge(\omega)\right|^p d\omega\right)^{1/p}$$

$$= \left(\int_{\mathbb{R}} |\omega|^{-pk} \left|\hat{\Psi}(\omega) - \left(\hat{\Psi}(\omega) - (\psi_n)^\wedge(\omega)\right)\right|^p d\omega\right)^{1/p}.$$

Using the Minkowski inequality, we get

$$\left(\int_{\mathbb{R}} |\omega|^{-pk} \left|\left((\psi_n)^\wedge(\omega) - \hat{\Psi}(\omega)\right) + \hat{\Psi}(\omega)\right|^p d\omega\right)^{1/p}$$

$$\leq \left(\int_{\mathbb{R}} |\omega|^{-pk} \left|\hat{\Psi}(\omega)\right|^p d\omega\right)^{1/p} + \left(\int_{\mathbb{R}} |\omega|^{-pk} \left|(\psi_n)^\wedge(\omega) - \hat{\Psi}(\omega)\right|^p d\omega\right)^{1/p} = I_1 + I_2.$$

On the other hand,

$$\left(\int_{\mathbb{R}} |\omega|^{-pk} \left|\left((\psi_n)^\wedge(\omega) - \hat{\Psi}(\omega)\right) + \hat{\Psi}(\omega)\right|^p d\omega\right)^{1/p}$$

$$\geq \left(\int_{\mathbb{R}} |\omega|^{-pk} \left|\hat{\Psi}(\omega)\right|^p d\omega\right)^{1/p} - \left(\int_{\mathbb{R}} |\omega|^{-pk} \left|(\psi_n)^\wedge(\omega) - \hat{\Psi}(\omega)\right|^p d\omega\right)^{1/p} = I_1 - I_2.$$

For $I_1^p$, we have

$$I_1^p = \left(\frac{1}{\sqrt{2\pi}}\right)^p \int_{-2\pi}^{-\pi} \omega^{-pk} d\omega + \left(\frac{1}{\sqrt{2\pi}}\right)^p \int_{\pi}^{2\pi} \omega^{-pk} d\omega = \frac{(2\pi)^{1/p-1/2}}{\pi^k} \left(\frac{1-2^{1-pk}}{pk-1}\right)^{1/p}.$$

Let us show that

$$I_2^p = \int_{\mathbb{R}} |\omega|^{-pk} \left|(\psi_n)^\wedge(\omega) - \hat{\Psi}(\omega)\right|^p d\omega \to 0 \quad \text{as} \quad n \to \infty.$$

We fix $\varepsilon \in (0; \pi)$. Dividing the interval of integration into two parts, we obtain

$$I_2^p = \int_{|\omega|<\varepsilon} |\omega|^{-pk} \left|(\psi_n)^\wedge(\omega) - \hat{\Psi}(\omega)\right|^p d\omega + \int_{|\omega|>\varepsilon} |\omega|^{-pk} \left|(\psi_n)^\wedge(\omega) - \hat{\Psi}(\omega)\right|^p d\omega = I_{11} + I_{12}.$$

Consider $I_{11}$. Taking into account condition (i) and the fact that $\hat{\Psi}(\omega) = 0$ for $\omega \in (-\pi; \pi)$, we establish that $I_{11} \to 0$ as $n \to \infty$. Furthermore, by virtue of condition (ii), we have

$$I_{12} \leq \varepsilon^{-kp} \int_{|\omega|>\varepsilon} |(\psi_n)^{\wedge}(\omega) - \hat{\Psi}(\omega)|^p d\omega \to 0 \quad \text{as} \quad n \to \infty.$$

Thus,

$$\|(_k\psi_n)^{\wedge}\|_p \to I_1^p = \frac{(2\pi)^{1/p-1/2}}{\pi^k}\left(\frac{1-2^{1-pk}}{pk-1}\right)^{1/p} \quad \text{as} \quad n \to \infty.$$

The lemma is proved.

We also note the special case of the lemma. For $k = 0$, we get

$$\lim_{n\to\infty} \|(\psi_n)^{\wedge}\|_q = (2\pi)^{1/q-1/2}. \tag{7}$$

**Proof of Theorem 1.** It is necessary to verify conditions (i) and (ii) of Lemma 1 for orthogonal Daubechies wavelets. We use relations (4) and (5) and the fact that

$$c_m = \left(\int_0^\pi \sin^{2m-1}\omega\, d\omega\right)^{-1} = \frac{\Gamma(m+1/2)}{\sqrt{\pi}\,\Gamma(m)} \sim \sqrt{\frac{m}{\pi}}. \tag{8}$$

To prove condition (i), we choose $0 < \varepsilon < 1$ and note that $|H_m(\omega)| \leq 1$ for any $\omega$. We obtain

$$\int_{|\omega|<\varepsilon} |\omega|^{-pk} |(\psi_m^D)^{\wedge}(\omega)|^p d\omega \leq \left(\frac{1}{2\pi}\right)^{p/2} \int_{|\omega|<\varepsilon} |\omega|^{-pk} \left|H_m\left(\frac{\omega}{2}+\pi\right)\right|^p d\omega$$

$$\leq \left(\frac{c_m}{2\pi}\right)^{p/2} \int_{|\omega|<\varepsilon} |\omega|^{-pk} \left(\int_0^{|\omega|/2} \sin^{2m-1}t\, dt\right)^{p/2} d\omega.$$

Furthermore, since

$$\left(\frac{|\omega|}{2}\right)^{-1} \sin\left(\frac{|\omega|}{2}\right) \leq 1,$$

we get

$$\left(\frac{c_m}{2\pi}\right)^{p/2} \int\limits_{|\omega|<\varepsilon} |\omega|^{-pk} \left(\int\limits_0^{|\omega|/2} \sin^{2m-1} t \, dt\right)^{p/2} d\omega$$

$$= 2^{-p(1/2+k)} \left(\frac{c_m}{\pi}\right)^{p/2} \int\limits_{|\omega|<\varepsilon} \left(\frac{|\omega|}{2}\right)^{-pk} \left(\int\limits_0^{|\omega|/2} \sin^{2m-1} t \, dt\right)^{p/2} d\omega$$

$$\leq 2^{-p(1/2+k)} \left(\frac{c_m}{\pi}\right)^{p/2} \int\limits_{|\omega|<\varepsilon} \left(\frac{|\omega|}{2}\right)^{-pk} \left(\frac{|\omega|}{2} \sin^{2m-1} \frac{|\omega|}{2}\right)^{p/2} d\omega$$

$$\leq 2^{-p(1/2+k)} \left(\frac{c_m}{\pi}\right)^{p/2} \int\limits_{|\omega|<\varepsilon} \left(\frac{|\omega|}{2}\right)^{-pk+p/2} \left(\sin^{2m-1} \frac{|\omega|}{2}\right)^{p/2} d\omega$$

$$\leq 2^{-p(1/2+k)} \left(\frac{c_m}{\pi}\right)^{p/2} \int\limits_{|\omega|<\varepsilon} \sin^{2mp/2-pk} \frac{|\omega|}{2} d\omega.$$

For $m > k$, we have

$$2^{-p(1/2+k)} \left(\frac{c_m}{\pi}\right)^{p/2} \int\limits_{|\omega|<\varepsilon} \sin^{2mp/2-pk} \frac{|\omega|}{2} d\omega \leq 2^{-p(1/2+k)} \left(\frac{c_m}{\pi}\right)^{p/2} 2\varepsilon \sin^{pm-pk} \frac{\varepsilon}{2}$$

$$\leq 2^{-p(1/2+k)+1} \left(\frac{c_m}{\pi}\right)^{p/2} \left(\frac{\varepsilon}{2}\right)^{pm-pk+1}.$$

Since

$$c_m \sim \sqrt{\frac{m}{\pi}},$$

the last expression tends to zero as $m \to \infty$.

Relation (i) is proved.

To prove condition (ii), we set $I = [-2\pi; 2\pi]$ and

$$I_\delta = [-2\pi, -2\pi + \delta) \cup (-\pi - \delta, -\pi + \delta) \cup (\pi - \delta, \pi + \delta) \cup (2\pi - \delta, 2\pi].$$

Let us prove that

$$\left(\int\limits_{I \setminus I_\delta} |(\psi_m^D)^\wedge(\omega) - \hat{\Psi}(\omega)|^p d\omega\right)^{1/p} \to 0 \quad \text{as} \quad m \to \infty.$$

We have

$$\left( \int_{I \setminus I_\delta} \left| (\psi_m^D)^\wedge(\omega) - \hat{\Psi}(\omega) \right|^p d\omega \right)^{1/p}$$

$$\leq \left( \int_{I \setminus I_\delta} \left| (\psi_m^D)^\wedge(\omega) - \frac{1}{\sqrt{2\pi}} H_m\left(\frac{\omega}{2} + \pi\right) \right|^p d\omega \right)^{1/p} + \left( \int_{I \setminus I_\delta} \left| \hat{\Psi}(\omega) - \frac{1}{\sqrt{2\pi}} H_m\left(\frac{\omega}{2} + \pi\right) \right|^p d\omega \right)^{1/p}. \quad (9)$$

For fixed $\delta$, the sequence $\frac{1}{\sqrt{2\pi}} H_m\left(\frac{\omega}{2} + \pi\right)$ converges uniformly as $m \to \infty$ to $\hat{\Psi}(\omega)$ in $I \setminus I_\delta$. Therefore, the second term on the right-hand side of (9) tends to zero as $m \to \infty$.

The first term on the right-hand side of (9) can be rewritten in the form

$$\frac{1}{2\pi} \left( \int_{I \setminus I_\delta} \left| H_m\left(\frac{\omega}{2} + \pi\right) \right|^p \left| 1 - \prod_{l \geq 1} \left| H_m(2^{-l-1}\omega) \right| \right|^p d\omega \right)^{1/p}. \quad (10)$$

Taking into account that the relation

$$\left| H_m\left(\frac{\omega}{2} + \pi\right) \right|^p \leq 1$$

and the estimates

$$\left| H_m\left(\frac{\omega}{4}\right) \right| \geq \left( 1 - c_m \frac{\pi}{2} \left\{ \sin\left(\frac{\pi}{2} - \frac{\delta}{4}\right) \right\}^{2m-1} \right)^{1/2},$$

$$\left| H_m\left(\frac{\omega}{8}\right) \right| \geq \left( 1 - c_m \left(\frac{\pi}{4}\right)^{2m} \right)^{1/2},$$

and

$$\prod_{l \geq 1} \left| H_m(2^{-l-3}\omega) \right|^p \geq \left( 1 - 2^{-2m} \right)^{\frac{1}{2(1-2^{-2m})}}$$

established in [4] for $\omega \in I \setminus I_\delta$, we can estimate integral (10) for all sufficiently large $m$ as follows:

$$\frac{1}{2\pi} \left( \int_{I \setminus I_\delta} \left| H_m\left(\frac{\omega}{2} + \pi\right) \right|^p \left| 1 - \prod_{l \geq 1} \left| H_m(2^{-l-1}\omega) \right| \right|^p d\omega \right)^{1/p}$$

$$\leq \frac{1}{2\pi} \left( (4\pi)^p \left| 1 - \left( 1 - c_m \frac{\pi}{2} \left\{ \sin\left(\frac{\pi}{2} - \frac{\delta}{4}\right) \right\}^{2m-1} \right)^{1/2} \left( 1 - c_m \left(\frac{\pi}{4}\right)^{2m} \right)^{1/2} \left( 1 - 2^{-2m} \right)^{\frac{1}{2(1-2^{-2m})}} \right| \right).$$

It is clear that the right-hand side of the inequality obtained tends to zero as $m \to \infty$.

Since $|I_\delta| = 6\delta,$ we have

$$\int_{I_\delta} |(\psi_m^D)^\wedge(\omega) - \hat{\Psi}(\omega)|^p d\omega \le \left(\frac{1}{2\pi}\right)^{p/2} |I_\delta| \le \left(\frac{1}{2\pi}\right)^{p/2} 6\delta. \tag{11}$$

We now verify that, for all $p > 1,$ we have

$$\int_{|\omega| \ge 2\pi} |(\psi_m^D)^\wedge(\omega)|^p d\omega \to 0 \quad \text{as} \quad m \to \infty.$$

According to known results concerning the regularity of Daubechies wavelets (see, e.g., [5], Sec. 2.2.4), there exist positive constants $C$ and $\tilde{C}$ such that the following inequality holds for all $\omega > 2\pi$:

$$|(\psi_m^D)^\wedge(\omega)| \le \tilde{C}|\omega|^{-C\log m}.$$

Then

$$\int_{|\omega|>2\pi} |(\psi_m^D)^\wedge(\omega)|^p d\omega \le \tilde{C} \int_{|\omega|>2\pi} |(\omega)|^{-Cp\log m} d\omega = \tilde{C} \int_{|\omega|>2\pi} |(\omega)|^{-(Cp\log m - 2)} d\omega$$

$$\le (2\pi)^{-(Cp\log m - 2)} \int_{|\omega|>2\pi} |(\omega)|^{-2} d\omega.$$

The last expression tends to zero as $m \to +\infty.$

Thus, the conditions of Lemma 1 for orthogonal Daubechies wavelets are satisfied.

Using (2), (6), and (7), we obtain

$$\lim_{m \to \infty} C_{k;p,q}(\psi_m^D) = \lim_{m \to \infty} \frac{\|(_k\psi_m^D)^\wedge\|_{p'}}{\|(\psi_m^D)^\wedge\|_q} = \frac{(2\pi)^{1/p - 1/2}}{\pi^k} \frac{\left(\frac{1 - 2^{1-pk}}{pk - 1}\right)^{1/p}}{(2\pi)^{1/q - 1/2}},$$

$$\lim_{m \to \infty} C_{k;p,q}(\psi_m^D) = \frac{(2\pi)^{1/p - 1/q}}{\pi^k} \left(\frac{1 - 2^{1-pk}}{pk - 1}\right)^{1/p}.$$

The theorem is proved.

## REFERENCES

1. I. Ya. Novikov and S. B. Stechkin, "Main theories of splashes," *Usp. Mat. Nauk,* **53**, No. 6, 53–128 (1998).
2. I. Daubechies, *Ten Lectures on Wavelets,* Society of Industrial and Applied Mathematics, Philadelphia (1992).
3. G. Strang and T. Nguyen, *Wavelets and Filter Banks,* Cambridge Press, Wellesley (1996).
4. S. Ehrich, "On the estimation of wavelet coefficients," *Adv. Comput. Math.,* **13**, 105–129 (2000).
5. A. K. Louis, P. Maab, and A. Rieder, *Wavelets: Theory and Applications,* Wiley, Chichester (1997).